\documentclass[a4paper,11pt]{article}
\usepackage{amsfonts}
\usepackage{graphics}
\newtheorem{theorem}{Theorem}[section]
\newtheorem{lemma}{Lemma}[section]
\newtheorem{definition}{Definition}[section]

\newtheorem{remark}{Remark}[section]
\def\dfrac{\displaystyle \frac}

\newcommand{\proof}{{\bf Proof:} }

\newcommand{\eop}{ \hfill $\Box$ }
\begin{document}
\begin{center}
{\Large Optimal Control of
Stochastic Functional Differential Equations with Application to Finance}
\\
\end{center}

\vspace{0.3cm}

\begin{center}
{\large  Edson A. Coayla-Teran}\footnote{Supported by CAPES Grant  4437/08-0} \\
\textit{Universidade Federal da Bahia-UFBA\\
Av. Ademar de Barros s/n, Instituto de Matem\' atica, Salvador, BA, Brasil, 40170-110, Telephone: 0055-71-32836279, e-mail: coayla@ufba.br}
\end{center}
\begin{center}
{\large Anatoly Swishchuk}\\
\textit{Department of
 Mathematics and Statistics, University of Calgary, Alberta, Canada,
 T2N1N4, Telephone: +1 (403) 220-3274 , e-mail: aswish@math.ucalgary.ca}
\end{center}
\begin{abstract}
This work is devoted to the study of optimal control of stochastic functional differential equations (SFDEs) and its application to mathematical finance. By using the Dynkin formula and solution of the Dirichlet-Poisson
problem, the Hamilton-Jacobi-Bellman (HJB) equation and the
converse HJB equation are derived. Furthermore, applications are given to
an optimal portfolio selection problem.
\end{abstract}
\section{Introduction} This work is devoted to the study of optimal control of stochastic functional differential equations (SFDEs). We believe that SFDEs are useful dynamical models to understand the behavior of natural process that take into consideration the influence of past events on the current and future states of the system \cite{SE2, IS, K}. This view is especially appropriate in the study of financial variables, since predictions about their evolution take strongly into account the knowledge of their past \cite{HO,SK}.

 The SFDEs are very important object that has many applications. One of the problems in the theory of SFDEs is the study of optimal control that has also many applications including finance. The main idea in finance is to find the optimal portfolio of an investor to maximize his wealth or cost
function. In this way, the SFDEs with controlled parameters are the main object of investigation of this paper.

The article is organized in the following way: in Section 2 we present the  basic spaces, the norms, properties and notation which we are going to work with in the following sections and formulation of the problem that is the goal of this work. In Section 3 we stated the results on existence and uniqueness of the solution of the SFDEs. We proved that the pair of processes, one with delayed parameter
and another one as the solution of the SFDEs, is a strong Markov process. With this result in hands and ussing  the weak infinitesimal generator of the Markov process (see \cite{YSE}, Lemma 9.3) we can apply the theory of controlled Markov processes to the solution of our optimization problem. We found the sufficient conditions for the optimality of the solution and derived the Hamilton-Jacobi-Bellman equation (HJB) equation and the converse of the HJB equation. In Section 4, the results obtained in Section 3 are applied to optimal portfolio selection problem where we found the optimal in explicit form.
\section{Preliminaries and formulation of the problem}
Let $a>0,$ $U$  be a closed set of $\mathbb{R}^m$ and $\left(\Omega,\mathcal{F},\left(\mathcal{F}_t\right)_{{0\leq t\leq a}},\mathbb{P}\right)$ be a complete filtered probability space. Assume also that each $\mathcal{F}_t$ contains all the sets of measure zero in $\mathcal{F}.$
Let  $r>0,$ $J:=[-r,0,]$  and $T:=[0,a].$\\
We denote by $V:=L^2\left([-r,0],\mathbb{R}^n\right),$ $H:=L^2\left([0,a],\mathbb{R}^n\right),$ with respective norms and inner products $\|\cdot\|_{V},$ $\left\langle\cdot,\cdot \right\rangle_{V},$ and $\|\cdot\|_{H},$ $\left\langle\cdot ,\cdot\right\rangle_{H}.$ Assume $\mu:V\times\mathbb{R}^n\times U\rightarrow \mathbb{R}^n,$ and $\sigma:V\times\mathbb{R}^n\times U\rightarrow \mathbb{R}^{n\times d}$ are measurable.
Now, we consider the following stochastic functional differential equation (SFDE)
\begin{equation}\label{SFDE}
S(t)=\left\{\begin{array}{rl}\stackrel{\rightarrow}{x}+\int_0^t&\!\!\!\!\mu(S_p,S(p),u(p))dp+\int_0^t\sigma(S_p,S(p),u(p))dW(p), t\in T\\
&\!\!\!\!\!\!\!\!\!\!\!\!\!\!\!\!\!\!\!\!\!\!\phi(t),\ \ -r\leq t< 0,
\end{array}\right.
\end{equation}
where $\phi$ is an initial path in $V,$ $\stackrel{\rightarrow}{x}$ an initial vector in $\mathbb{R}^n$ and $W(t)$ is an  $\mathcal{F}_t-$ adapted $d-$dimensional Brownian motion, and $\mathbf{u}(s)$ is defined below.\\
The solution $\left\{S(t)\right\}_{-r\leq t\leq a}$ of (\ref{SFDE}) is an $n-$dimensional stochastic process. Its \textit{segment process} $\left\{S_t: t\in T\right\}$ is defined by
\begin{equation}\label{segpr}
S_t(\omega)(p):=S(t+p,\omega) \textmd{ for }p\in J.
\end{equation}
The function $u(t):=\mathbf{u}(t)=\mathbf{u}(t,S_t, S(t))$ will be called \textit{Markov control law}. A Markov control law $\mathbf{u}:J\times\times\mathbb{R}^n\rightarrow U$ is \textit{admissible} if it is a Borel  measurable function and it satisfies:
$$|\mathbf{u}(t,\phi,\stackrel{\rightarrow}{x})-\mathbf{u}(t,\eta,\stackrel{\rightarrow}{y})|^2\leq K\left\{|\stackrel{\rightarrow}{x}-\stackrel{\rightarrow}{y}|^2+\|\phi-\eta\|_V^2\right\}$$  for some constant $K,$ for all $\eta,$ $\phi\in V,$ $t\in J,$ $\stackrel{\rightarrow}{x},\stackrel{\rightarrow}{y}\in\mathbb{R}^n$  and holds
$$|\mathbf{u}(t,\phi,\stackrel{\rightarrow}{x})|^2\leq K_1\left\{1+|\stackrel{\rightarrow}{x}|^2+\|\phi\|_V^2\right\},$$ for some constant $K_1$ for all $\phi\in V,$ $t\in J,$ $\stackrel{\rightarrow}{x}\in\mathbb{R}^n.$ We denote by $\mathcal{U}$ the set of all admissible Markov control laws.\\
Let $G\subseteq V\times\mathbb{R}^n$ be a open connected subset with boundary $\Gamma:=\partial(G).$  Let $\psi(\cdot,\cdot)$ be a continuous function  on the closure of the set $G$ and bounded on $\Gamma$ and  $L(\cdot,\cdot,\cdot)$ be a continuous function on $G\times U$ be such that\\
\begin{equation}\label{condint}
E^{(\phi,\stackrel{\rightarrow}{x})}\left[\int_0^{\tau_{G}}|L(S_t,S(t),u(t))|dt\right]<\infty\ \ \forall\,(\phi,\stackrel{\rightarrow}{x})\in G,
\end{equation}
where $\tau_{G}$ is the first exit time from the set $G,$ and $E^{({\phi,\stackrel{\rightarrow}{x}})}$ the expectation with respect the probability laws $\mathbb{P}^{({\phi,\stackrel{\rightarrow}{x}})}$ of $(S_t,S(t)),$ for $(\phi,\stackrel{\rightarrow}{x})\in G$.\\
Now we are given a \textit{cost function}(or \textit{performance criterion})
\begin{equation}\label{PC}\begin{array}{rl}
J(\phi,\stackrel{\rightarrow}{x},\mathbf{u}):=&\!\!\!E^{(\phi,\stackrel{\rightarrow}{x})}\left[\int_0^{\tau_{G}}L(S_t,S(t),u(t))dt\right.+\\
&\!\!\!+\left.\psi(S_{\tau_{G}},S(\tau_{G}))\right].
\end{array}
\end{equation}
The problem is to find the number $\Phi(\phi,\stackrel{\rightarrow}{x})$ and a control $\mathbf{u}^{\star}=\mathbf{u}^{\star}(t,\omega),$ for each $\left(\phi,\stackrel{\rightarrow}{x}\right)\in G,$   such that
\begin{equation}\label{infprob}
\Phi(\phi,\stackrel{\rightarrow}{x}):=\inf_{\mathbf{u}}J(\phi,\stackrel{\rightarrow}{x},\mathbf{u})=J(\phi,\stackrel{\rightarrow}{x},\mathbf{u}^{\star})
\end{equation}
where the infimum is taken over all $\mathcal{F}_t-$adapted process $\mathbf{u}\in\mathcal{U}$.  If a such control $\mathbf{u}^{\star}$  exists then it is called an \textit{optimal control} and $\Phi$ is called the \textit{optimal performance}.\\
We denote by $B_b(V\times \mathbb{R}^n)$ the Banach space of all real bounded Borel functions, endowed with the sup norm.
\section{Controlled Stochastic Differential Delay Equations}
Given the Markov control $u(t)=\mathbf{u}(S_t,S(t))$ and a function $g(\phi,\stackrel{\rightarrow}{x},u),$ we use the notation
$$
g^{\mathbf{u}}\left(\phi,\stackrel{\rightarrow}{x}\right)=g(\phi,\stackrel{\rightarrow}{x},\mathbf{u}(\phi,\stackrel{\rightarrow}{x})).
$$
Then (\ref{SFDE}) can be write as
\begin{equation}\label{SFDEu}
S(t)=\left\{\begin{array}{rl}\stackrel{\rightarrow}{x}+\int_0^t&\!\!\!\!\mu^{\mathbf{u}}(S_p,S(p))dp+\int_0^t\sigma^{\mathbf{u}}(S_p,S(p))dW(p), t\in T\\
&\!\!\!\!\!\!\!\!\!\!\!\!\!\!\!\!\!\!\!\!\!\phi(t),\ \ -r\leq t< 0,
\end{array}\right.
\end{equation}
\begin{theorem}\label{lexuq}Let $\phi:\Omega\rightarrow V$ such that $E\left[\|\phi\|_V^2\right]<+\infty$ and $\stackrel{\rightarrow}{x}:\Omega\rightarrow \mathbb{R}^n$ such that $E\left[\|\stackrel{\rightarrow}{x}\|^2\right]<+\infty$ and $\mathcal{F}_0$ mensurable. Assume that there exists a constant $K$ such that
\begin{equation}\label{exstcond1}\begin{array}{rl}
\left\|\mu^{\mathbf{u}}(\phi,\stackrel{\rightarrow}{x})\right.&\!\!\!-\left.\mu^{\mathbf{u}}(\eta,\stackrel{\rightarrow}{x_1})\right\|^2+\left\|\sigma^{\mathbf{u}}(\phi,\stackrel{\rightarrow}{x})-\sigma^{\mathbf{u}}(\eta,\stackrel{\rightarrow}{x_1})\right\|^2\leq\\
&\leq K\left[\left\|\stackrel{\rightarrow}{x}-\stackrel{\rightarrow}{x_1}\right\|^2+\|\phi-\eta\|_V^2\right]
\end{array}
\end{equation}
and
\begin{equation}\label{exstcond2}
|\mu^{\mathbf{u}}(\phi,\stackrel{\rightarrow}{x},u)|^2+|\sigma^{\mathbf{u}}(\phi,\stackrel{\rightarrow}{x},u)|^2\leq K(1+|\stackrel{\rightarrow}{x}|^2+\|\phi\|^2_V).
\end{equation}
for all $\phi,\ \eta\in V,$ $\stackrel{\rightarrow}{x},\ \stackrel{\rightarrow}{x_1}\in\mathbb{R}^n.$ \\
Then we have a unique measurable solution $S(t)$ to (\ref{SFDEu}) with continuous trajectories $\left\{(S_t,S(t)),\ t\in T\right\}$ adapted to $\left(\mathcal{F}_t\right)_{t\in T}$.
\end{theorem}
\proof The proof is by using the method of successive approximations (see \cite{SE3}, page 227).\eop\\
To the case $n=1$ we can still assure the existence and uniqueness of solution to (\ref{SFDEu}) under weaker conditions.
\begin{theorem}\label{lexuq1}Under the same notations of Theorem (\ref{lexuq}) and $n=1$. Let $\phi:\Omega\rightarrow V$ such that $E\left[\|\phi\|_V^2\right]<+\infty$ and $\stackrel{\rightarrow}{x}:\Omega\rightarrow \mathbb{R}$ such that $E\left[|\stackrel{\rightarrow}{x}|^2\right]<+\infty$ and $\mathcal{F}_0-$mensurable. Assume that there exists a constant $K$ such that
\begin{equation}\label{exstcond2}
|\mu^{\mathbf{u}}(\phi,\stackrel{\rightarrow}{x})|^2+|\sigma^{\mathbf{u}}(\phi,\stackrel{\rightarrow}{x})|^2\leq K^2(1+|\stackrel{\rightarrow}{x}|^2+\|\phi\|^2_V).
\end{equation}
for all $\phi\in V,$ $\stackrel{\rightarrow}{x}\in\mathbb{R},$ $u\in U.$\\
And for each $N$ there exists $K_N$ for which
\begin{equation}\label{exstcond1}\begin{array}{rl}
\left|\mu^{\mathbf{u}}(\phi,\stackrel{\rightarrow}{x})\right.&\!\!\!-\left.\mu^{\mathbf{u}}(\eta,\stackrel{\rightarrow}{x_1})\right|^2+\left|\sigma^{\mathbf{u}}(\phi,\stackrel{\rightarrow}{x})-\sigma^{\mathbf{u}}(\eta,\stackrel{\rightarrow}{x_1})\right|^2\leq\\
&\leq K_N\left[\left|\stackrel{\rightarrow}{x}-\stackrel{\rightarrow}{x_1}\right|^2+\|\phi-\eta\|_V^2\right]
\end{array}
\end{equation}
for all $\phi,\ \eta\in V,$ $\stackrel{\rightarrow}{x},\ \stackrel{\rightarrow}{x_1}\in\mathbb{R},$ with $|\stackrel{\rightarrow}{x}|\leq K_N,$ $|\stackrel{\rightarrow}{x_1}|\leq K_N$\\
Then we have a unique measurable solution $S(t)$ to (\ref{SFDEu}) with continuous trajectories $\left\{(S_t,S(t)),\ t\in T\right\}$ adapted to $\left(\mathcal{F}_t\right)_{t\in T}$.
\end{theorem}
\proof See \cite{GS}, Theorem 3, page 45.\eop\\
\begin{remark}\label{exuqarb}Let $0\leq t_1\leq t\leq T,$  $\phi\in V,$  $\stackrel{\rightarrow}{x}:\Omega\rightarrow \mathbb{R}^n$ such that $E\left[\|\stackrel{\rightarrow}{x}\|^2\right]<+\infty$ with $\phi,\  \stackrel{\rightarrow}{x}\mathcal{F}_{t_1}-$mensurable. We can solve the following equation at time $t_1$
\begin{equation}
\left\{\begin{array}{rl}\label{SDDEs}
S(t)&\!\!\!=\stackrel{\rightarrow}{x}+\int_{t_1}^t\mu^{\mathbf{u}}(S_p,S(p))dp+\int_{t_1}^t\sigma^{\mathbf{u}}(S_p,S(p))dW(p), t\in[ t_1,T]\\
S(t)&\!\!\!=\phi(t-t_1),t\in [t_1-r,t_1).\end{array}\right.
\end{equation}
We denote by $S(\cdot,t_1,\phi,\stackrel{\rightarrow}{x})$ the solution of (\ref{SDDEs}).
\end{remark}
Moreover, the solution have similar properties that the solutions of stochastic differential equations.
\begin{theorem}\label{tpropr}Under the assumptions of Theorem \ref{lexuq}  , there exists $C(a,r)>0$ such that, for arbitrary $\phi,\eta:\Omega\rightarrow V$ such that $E\left[\|\phi\|_V^2\right],\ E\left[\|\eta\|_V^2\right]<+\infty,$ and $x,y:\Omega\rightarrow \mathbb{R}^n$ such that $E\left[\|\stackrel{\rightarrow}{x}\|^2\right],\ E\left[\|\stackrel{\rightarrow}{y}\|^2\right]<+\infty$ and $\mathcal{F}_0$ mensurable. Let $0\leq t_1\leq t\leq a,$ then
\begin{equation}
E\left(\|(S(\cdot,t_1,\phi,\stackrel{\rightarrow}{x}))_t\|_{V}^2+|S(t,t_1,\phi,\stackrel{\rightarrow}{x})|^2\right)\leq C(a,r)E(\|\phi\|^2_V)+E\left(|\stackrel{\rightarrow}{x}|^2\right)
\end{equation}
\begin{equation}\begin{array}{rl}
\sup_{s\in[t_1,a]}E\left(|S(s,\right.&\!\!\!\!\!\left.t_1,\phi,\stackrel{\rightarrow}{x})-(S(s,t_1,\eta,\stackrel{\rightarrow}{y}))|\right)+\\
&\!\!+E\left(\|(S(\cdot,t_1,\phi,\stackrel{\rightarrow}{x}))_t-(S(\cdot,t_1,\eta,\stackrel{\rightarrow}{y}))_t\|_{V}^2\right)\leq\\
&C(a,r)E\left(\|\phi-\eta\|_V^2\right)+E\left(|\stackrel{\rightarrow}{x}-\stackrel{\rightarrow}{y}|^2\right).
\end{array}
\end{equation}
\begin{equation}\begin{array}{rl}
E\left(|S(t,t_1,\phi,\stackrel{\rightarrow}{x})-\right.&\!\!\!\!\!\left.S(t,t_1,\phi,\stackrel{\rightarrow}{x})|\right)+\\
&\!\!\!+E\left(\|(S(\cdot,t_1,\phi,\stackrel{\rightarrow}{x}))_t-(S(\cdot,t_1,\eta,\stackrel{\rightarrow}{y}))_t\|_{V}^2\right)\leq\\
&\!\!\!C(a,r,\phi,\stackrel{\rightarrow}{x})|t-t_1|^2.
\end{array}
\end{equation}
\end{theorem}
\proof The proof is using similar ideas as in the case of no delay (see \cite{DZ}, Theorem 9.1) and similat to Theorem 3.1, page 41 from \cite{SE3}.\eop\\
Let $A\in \mathcal{B}(\mathbb{R}^n)\otimes\mathcal{B}(V),$ we define the transition probability
$$\begin{array}{rl}
p\left(t_1,(\phi,\stackrel{\rightarrow}{x}),t,A\right):=&\mathbb{P}\left(\left((S(\cdot,t_1,\phi,\stackrel{\rightarrow}{x}))_t,S(t,t_1,\phi,\stackrel{\rightarrow}{x})\right)\in A\right)=\\
&=E\left[1_A(S(\cdot,t_1,\phi,\stackrel{\rightarrow}{x})_t,S(t,t_1,\phi,\stackrel{\rightarrow}{x}))\right].
\end{array}
$$
We will show now, following \cite{GS}, that the process $\left(S_t,S(t)\right),$ $t\in T$, is a Markov process with transition probability $p\left(t_1,(\phi,\stackrel{\rightarrow}{x}),t,A\right).$
\begin{lemma}\label{markv} Assume that $S(t)$ is solution to (\ref{SFDEu}). Then $\left(S_t,S(t)\right)$ will be a Markov process with transition probability $p\left(t_1,(\phi,\stackrel{\rightarrow}{x}),t,A\right),$ $0\leq t_1\leq t\leq a.$ and $A\in \mathcal{B}(\mathbb{R}^n)\otimes\mathcal{B}(V).$
\end{lemma}
\proof
Denote by $G^t$ the $\sigma-$algebra generated by $W(s)-W(t)$ for $t\leq s.$ We observe that $G^t$ and $\mathcal{F}_t$ are independent.
We observe that $S(t)= S(t,t_1,S_{t_1},S(t_1))$ for $t>t_1$, because both are solutions of the equation:
\begin{equation}\left\{\begin{array}{rl}
Z(t)&\!\!\!=Z(t_1)+\int_{t_1}^t\mu^{\mathbf{u}}(Z_s,Z(s))ds+\int_{t_1}^t\sigma^{\mathbf{u}}(Z_s,Z(s))dW(s),\ t_1\leq t\leq T\\
Z(t)&\!\!\!=S(t-t_1),t\in [t_1-r,t_1)
\end{array}\right.
\end{equation}
Let $B\in \mathcal{F}_{t_1}.$ Since
$$\begin{array}{rl}
\int_{B}1_{A}((&\!\!\!\!S_t,S(t)))d\mathbb{P}(\omega)=\int_{\Omega}1_{A}((S_t,S(t)))1_Bd\mathbb{P}(\omega)=\\
&\!\!\!=\int_{\Omega}1_{A}((S(\cdot,t_1,S_{t_1},S(t_1))_t,S(t,t_1,S_{t_1},S(t_1)))1_Bd\mathbb{P}(\omega)=\\
&\!\!\!=\int_{\Omega}1_{A}((S(\cdot,t_1,S_{t_1},S(t_1))_t,S(t,t_1,S_{t_1},S(t_1))))d\mathbb{P}(\omega)\int_{\Omega}1_Bd\mathbb{P}(\omega)=\\
&\!\!\!=\int_{B}\mathbb{P}\left(t_1,(\phi,\stackrel{\rightarrow}{x}),t,A\right)d\mathbb{P}(\omega)|_{x=S(t_1),\phi=S_{t_1}},\\
\end{array}
$$
thus we have that $\mathbb{P}\left((S_t,S(t))\in A|\mathcal{F}_{t_1}\right)=p(t_1,(S_{t_1},S{t_1}),t,A).$
To see that $\mathbb{P}\left((S_t,S(t))\in A|(S_{t_1},S(t_1))\right)=p(t_1,(S_{t_1},S(t_1)),t,A)$, we prove first that $P(t_1,.,t,A)$ is measurable for fixed $t,t_1,A,$ since $(S_{t_1},S(t_1))$ is measurable with respect to $\sigma-$algebra generated by $(S_{t_1},S(t_1))$
we finish the proof.
\eop\\
With similar arguments we can prove the following theorem. See for example \cite{DZ}, Theorem 9.8.
\begin{theorem}\label{markprop}Let  $S(t):=S(t,t_1,\phi,\stackrel{\rightarrow}{x})$ be the solution to (\ref{SDDEs}). For arbitrary $f\in B_b(V\times \mathbb{R}^n)$ and $0\leq t_1\leq t\leq a,$
\begin{equation}E\left[f(S_t,S(t))|\mathcal{F}_{t_1}\right]=E\left[f(S_t,S(t))|(S_{t_1},S(t_1))\right]
\end{equation}
\end{theorem}
\eop\\
Now, following \cite{DZ} we will prove that the solutions to (\ref{SFDEu}) are a strong Markov process.
\begin{theorem}\label{strgmrkv}(The strong Markov property) Let $S(t)$  as in the Theorem \ref{markprop},  $f$ in $B_b\left(V\times \mathbb{R}^n\right),$ $\tau$ a stopping time with respect to $\mathcal{F}_t,$ $\tau<\infty$ $\textit{a.s.}$ Then
\begin{equation}\label{strgmrq}
E^{(\phi,\stackrel{\rightarrow}{x})}\left[f(S_{\tau+h},S(\tau+h))|\mathcal{F}_{\tau}\right]=E^{(S_{\tau},S(\tau))}f(S_{h},S(h))
\end{equation}
for all $h\geq0.$
\end{theorem}
\proof We prove (\ref{strgmrq}) as in \cite{DZ}, Theorem 9.14 page 255 using the properties of Theorem \ref{tpropr}.
\eop\\
For every $f\in B_b\left(V\times\mathbb{R}^n\right)$ and $(\phi,\stackrel{\rightarrow}{x})\in V\times\mathbb{R}^n$ let
$$P_tf(\phi,\stackrel{\rightarrow}{x}):=E^{(\phi,\stackrel{\rightarrow}{x})}\left(f(S_t,S(t))\right).$$
\begin{definition}The weak infinitesimal operator of $P_t$ (or of $(S_t,S(t))$), $\mathbf{A}^{\mathbf{u}}:=\mathbf{A}^{\mathbf{u}}_S,$  is defined by
\begin{equation}\label{inft}
\mathbf{A}^{\mathbf{u}}f(\phi,\stackrel{\rightarrow}{x}):=\lim_{h\rightarrow0}h^{-1}\left[P_{h}f(\phi,\stackrel{\rightarrow}{x})-f(\phi,\stackrel{\rightarrow}{x})\right].
\end{equation}
The set  of functions $f$ such that the limit (\ref{inft}) exists in $(\phi,\stackrel{\rightarrow}{x})$ is denoted by  $\mathcal{D}_{\mathbf{A}^{\mathbf{u}}}(\phi,\stackrel{\rightarrow}{x})$ and $\mathcal{D}_{\mathbf{A}^{\mathbf{u}}}$ denotes the set of functions  such that the limit (\ref{inft}) exists for all $(\phi,\stackrel{\rightarrow}{x}).$
\end{definition}
Let $e_j$  for $j=1,\ldots,d$ be the canonical basis of $\mathbb{R}^d$ for $(\phi,\stackrel{\rightarrow}{x})\in V\times\mathbb{R}^n$ let
\begin{equation}\label{defS}
\widehat{\phi^{\stackrel{\rightarrow}{x}}}(t):=\left\{\begin{array}{rl}
\stackrel{\rightarrow}{x},&t\in T\\
\phi(t),&t\in [-r,0)
\end{array}\right.
\end{equation}
Then, for each $s\in J,$ $t\in T,$ 
\begin{equation}\widehat{\phi^{\stackrel{\rightarrow}{x}}_t}(s)=\widehat{\phi^{\stackrel{\rightarrow}{x}}}(s+t)=\left\{\begin{array}{rl}
\stackrel{\rightarrow}{x},&t+s\geq0\\
\phi(t),&t+s<0
\end{array}\right.
\end{equation}\\
Denote by $\mathbf{\Gamma}_t$ for $t\in T$ the weakly continuous contraction semigroup of the shift operators defined on $C_b(V\times\mathbb{R}^n)$ (see \cite{SE3}, Chapter 4) by
$$
\mathbf{\Gamma}_t(f)(\phi,\stackrel{\rightarrow}{x}):=f(\widehat{\phi^{\stackrel{\rightarrow}{x}}_t},\stackrel{\rightarrow}{x}) \textmd{ for } f\in C_b(V\times\mathbb{R}^n)
$$
Denote by $\mathbf{\Gamma}$ the weak infinitesimal operator of $\mathbf{\Gamma}_t$ with domain $D(\mathbf{S})$ and $D(\mathbf{S})\subset C^0_b=\left\{f\in C_b(V\times\mathbb{R}^n): S_t \textmd{ is strongly continuous}\right\}.$   
Now we have a formula for the weak infinitesimal operator $\mathbf{A}^{\mathbf{u}}$ similar to no delay case this is a sum of differential operators and depend of the coefficients $\mu^{\mathbf{u}}$ and $\sigma^{\mathbf{u}}.$
\begin{theorem}\label{wkinfopr} Let $S(t)$ be the solution to (\ref{SFDEu}). Suppose $f\in C^2_b(V\times\mathbb{R}^n)$, belongs to the domain of $\mathbf{A}^{\mathbf{u}},$  $\sigma^i\in C^2_b(V\times\mathbb{R}^n\times U;\mathbb{R}^n)$ (where $\sigma^i$ are the vector columns of $\sigma$) and $\mu\in C^1_b(V\times\mathbb{R}^n\times U;\mathbb{R}^n).$ Assume that $\phi\in V,$ $\stackrel{\rightarrow}{x}\in\mathbb{R}^n.$ Let ${e_j:j=1,\ldots,d}$ be a normalized basis of $\mathbb{R}^d.$ Then
\begin{equation}\label{wkinfop}\begin{array}{rl}
\mathbf{A}^{\mathbf{u}}f(\phi,\stackrel{\rightarrow}{x})=&\!\!\!\!\mathbf{\Gamma}f(\phi,\stackrel{\rightarrow}{x})+\frac{\partial f}{\partial \stackrel{\rightarrow}{x}}(\phi,\stackrel{\rightarrow}{x})\mu^{\mathbf{u}}(\phi,\stackrel{\rightarrow}{x})+\\
&\!\!\!+\frac{1}{2}\sum_j^n\frac{\partial^2 f}{\partial \stackrel{\rightarrow}{x}^2}(\phi,\stackrel{\rightarrow}{x})\left[(\sigma^{\mathbf{u}}(\phi,\stackrel{\rightarrow}{x}))e_j\otimes(\sigma^{\mathbf{u}}(\phi,\stackrel{\rightarrow}{x}))e_j\right]\\
\end{array}
\end{equation}
\end{theorem}
\proof Is consequence of Lemma 9.3 of \cite{YSE}.
\eop\\
\begin{remark}\label{operL}
Let $\mathbb{L}$ denote the differential operator given by the right hand side of (\ref{wkinfop}). The Theorem \ref{wkinfopr} above says that  $\mathbf{A}^{\mathbf{u}}$ and $\mathbb{L}$ coincide on $f\in C^2_b(V\times\mathbb{R}^n)$.
\end{remark}
\begin{lemma}(Dynkin formula). \label{dynkform}Let $S(t)$ be the solution of (\ref{SFDEu}). Let $f\in C^2_b(V\times\mathbb{R}^n),$ $\tau$ is a stopping time such that $E^{(\phi,\stackrel{\rightarrow}{x})}\left[\tau\right]<\infty,$ with $(\phi,\stackrel{\rightarrow}{x})\in V\times\mathbb{R}^n$ then
\begin{equation}
E^{(\phi,\stackrel{\rightarrow}{x})}\left[f(S_\tau,S(\tau))\right]=f(\phi,\stackrel{\rightarrow}{x})+E^{(\phi,\stackrel{\rightarrow}{x})}\left[\int_0^\tau\mathbf{A}^{\mathbf{u}}f(S_s,S(s))ds\right]
\end{equation}
\end{lemma}
\proof From Dynkin \cite{D}, corollary of Theorem 5.1.\eop\\
\begin{definition}Let $S(t)$ the solution of (\ref{SFDEu}). The characteristic operator $\mathcal{A}^{\mathbf{u}}=\mathcal{A}^{\mathbf{u}}_S$ of $(S_t,S(t))$ is defined by
\begin{equation}\label{charop}
\mathcal{A}^{\mathbf{u}}f(\phi,\stackrel{\rightarrow}{x}):=\lim_{U\downarrow (\phi,\stackrel{\rightarrow}{x})}\frac{E^{(\phi,\stackrel{\rightarrow}{x})}\left[f(S_{\tau_U},S(\tau_U))\right]-f(\phi,\stackrel{\rightarrow}{x})}{E^{\phi,\stackrel{\rightarrow}{x}}\left[\tau_U\right]}
\end{equation}
where the $U^{'}$s are open sets $U_k$ decreasing to the point $(\phi,\stackrel{\rightarrow}{x}),$ in the sense that $U_{k+1}\subset U_k$ and $\bigcap_k U_k={(\phi,\stackrel{\rightarrow}{x})},$ and $\tau_U=\inf\left\{t>0;\ (S_t,S(t))\notin U\right\}.$ We denote by $\mathcal{D}_{\mathcal{A}^{\mathbf{u}}}$ the set of functions $f$ such that the limit (\ref{charop}) exists for all $(\phi,\stackrel{\rightarrow}{x})\in V\times\mathbb{R}^n$ (and all $\left\{U_k\right\}$.) If $E^{(\phi,\stackrel{\rightarrow}{x})}\left[\tau_U\right]=\infty$ for all open $U\ni(\phi,\stackrel{\rightarrow}{x}),$ we define $\mathcal{A}^{\mathbf{u}}f(\phi,\stackrel{\rightarrow}{x})=0.$
\end{definition}
\begin{theorem}
Let $f\in C^2(V\times\mathbb{R}^n).$ Then $f\in\mathcal{D}_{\mathcal{A}^{\mathbf{u}}}$ and
\begin{equation}\label{coincAL}
\mathcal{A}^{\mathbf{u}}f=\mathbb{L}f.
\end{equation}
Where $\mathbb{L}$ is defined in Remark \ref{operL}.
\end{theorem}
\proof See \cite{O}, Theorem 7.5.4.
\eop\\
\begin{theorem}\label{possdrc}Assume that $\tau_{G}<\infty$ a.s. $\mathbb{P}^{(\phi,\stackrel{\rightarrow}{x})}.$ for all ${(\phi,\stackrel{\rightarrow}{x})}\in G$. Let $\psi\in C(\partial (G))$ be  bounded and let $g\in C(G)$ satisfy
\begin{equation}\label{cndpossdrh}
E^{(\phi,\stackrel{\rightarrow}{x})}\left[\int_0^{\tau_{G}}|g(S_t,S(t))|dt\right]<\infty,\ \forall\ (\phi,\stackrel{\rightarrow}{x})\in G.
\end{equation}
Define
\begin{equation}\begin{array}{rl}
w(\phi,\stackrel{\rightarrow}{x})&=E^{(\phi,\stackrel{\rightarrow}{x})}\left[\psi(S_{\tau_{G}},S(\tau_{G}))\right]+\\
&+E^{(\phi,\stackrel{\rightarrow}{x})}\left[\int_0^{\tau_{G}}g(S_t,S(t))dt\right],\ (\phi,\stackrel{\rightarrow}{x})\in G.
\end{array}
\end{equation}
Then\\
a) \begin{equation}\label{eqposs}
\mathcal{A}^{\mathbf{u}}w=-g\textmd{ in }\, G
\end{equation}
and
\begin{equation}\label{posscond}
\lim_{t\uparrow\tau_{G}}w(S_t,S(t))=\psi(S_{\tau_{G}},S(\tau_{G}))\ a.s.,
\end{equation}
b) Moreover, if there exists a function $w_1\in C^2(G)$ and a constant $C$ such that
\begin{equation}
|w_1(\phi,\stackrel{\rightarrow}{x})|<C\left(1+E^{(\phi,\stackrel{\rightarrow}{x})}\left[\int_0^{\tau_{G}}|g(S_t,S(t))|dt\right]\right),\ \left(\phi,\stackrel{\rightarrow}{x}\right)\in G,
\end{equation}
and $w_1$ satisfies (\ref{eqposs}) and (\ref{posscond}), then $w_1=w.$
\end{theorem}
\proof The proof follows similar arguments as \cite{O} Theorem 9.3.3.
\eop\\
\indent Let $M:V\times\mathbb{R}^n\times U\rightarrow \mathbb{R},$ such that $E^{(\phi,\stackrel{\rightarrow}{x})}\int_0^{\tau_{G}}|M^{\mathbf{u}}(S_t,S(t))|dt<\infty,$ we consider the equation
\begin{equation}\label{edpop}
(\mathbf{A}^{\mathbf{u}}f+M^{\mathbf{u}})(\phi,\stackrel{\rightarrow}{x})=0,\ \ (\phi,\stackrel{\rightarrow}{x})\in G
\end{equation}
with boundary data
\begin{equation}\label{boun}
f(\phi,\stackrel{\rightarrow}{x})=\psi(\phi,\stackrel{\rightarrow}{x})\textmd{ with } (\phi,\stackrel{\rightarrow}{x})\in\partial^{\ast}(G).
\end{equation}
Here $\partial^{\ast}(G) $ denotes a closed subset of $\partial (G)$ such that $\mathbb{P}^{(\phi,\stackrel{\rightarrow}{x})}((S_{\tau_{G}},S{\tau_{G}})\notin\partial^{\ast}(G),\ \tau_{G}<\infty)=0$ for each $(\phi,\stackrel{\rightarrow}{x})\in G.$
\begin{lemma}\label{lemut}Let $S(t)$ be the solution to (\ref{SFDEu}), $f$ in $C^2(G),$ with $F$ continuous and bounded. Suppose that $\mathbb{P}^{(\phi,\stackrel{\rightarrow}{x})}\left(\tau_{G}<\infty\right)=1$ for each $(\phi,\stackrel{\rightarrow}{x})\in G.$\\
\indent(a) If ($\mathbf{A}^{\mathbf{u}}f+M^{\mathbf{u}})(\phi,\stackrel{\rightarrow}{x})\geq0$ for all $(\phi,\stackrel{\rightarrow}{x})\in G,$ then
\begin{equation}\label{lemop}
f(\phi,\stackrel{\rightarrow}{x})\leq E^{(\phi,\stackrel{\rightarrow}{x})}\left\{\int_0^{\tau_{G}}M^{\mathbf{u}}(S_t,S(t))dt+f(S_{\tau_{G }},S(\tau_{G})\right\},\ \ (\phi,\stackrel{\rightarrow}{x})\in G
\end{equation}
\indent (b) If $f$ is a solution of (\ref{edpop}) and (\ref{boun}) for all $(\phi,\stackrel{\rightarrow}{x})\in G,$ where $E^{(\phi,\stackrel{\rightarrow}{x})}\int_0^{\tau_{G}}|M^{\mathbf{u}}(S_t,S(t))|<\infty,$ then
\begin{equation}\label{lemop}
f(\phi,\stackrel{\rightarrow}{x})= E^{(\phi,\stackrel{\rightarrow}{x})}\left\{\int_0^{\tau_{G}}M^{\mathbf{u}}(S_t,S(t))dt+\Psi(S_{\tau_{G}},S(\tau_{G})\right\},\ \ (\phi,\stackrel{\rightarrow}{x})\in G
\end{equation}
\end{lemma}
\proof (a) From Dynkin formula
$$\begin{array}{rl}
f(\phi,\stackrel{\rightarrow}{x})=&E^{(\phi,\stackrel{\rightarrow}{x})}f(S_{\tau_{G}},S(\tau_{G})+\\
&-E^{(\phi,\stackrel{\rightarrow}{x})}\left\{\int_0^{\tau_{G}}\mathbf{A}^{\mathbf{u}}f(S_t,S(t))dt\right\}\leq\\
&\leq E^{(\phi,\stackrel{\rightarrow}{x})}\left\{\int_0^{\tau_{G}}M^{\mathbf{u}}(S_t,S(t))dt+f(S_{\tau_{G}},S(\tau_{G})\right\}
\end{array}
$$
(b) Since $M^{\mathbf{u}}=-\mathbf{A}^{\mathbf{u}}f$ satisfies the condition integrability, we get (b) as in (a).
\eop\\
For $v=\mathbf{u}(S_t,S(t)),$ let
$$
\mathbf{A}^{v}f(S_t,S(t)):=\mathbf{A}^{\mathbf{u}}f(S_t,S(t))
$$
The dynamic programming equation is:
\begin{equation}\label{eqprg}
0=\inf_{v\in U}\left[(\mathbf{A}^{v}f+L^v)(\phi,\stackrel{\rightarrow}{x})\right],\ (\phi,\stackrel{\rightarrow}{x})\textmd{ in }G,
\end{equation}
with the boundary data
\begin{equation}\label{boun1}
f(\phi,\stackrel{\rightarrow}{x})=\psi(\phi,\stackrel{\rightarrow}{x})\ \ (\phi,\stackrel{\rightarrow}{x})\in\partial^{\ast}(G),
\end{equation}
and $L$ as in (\ref{condint}).\\
We assume that
\begin{equation}\label{inflim}
L(\phi,\stackrel{\rightarrow}{x},v)\geq c>0
\end{equation}
for some constant $c.$\\
One of the fundamental results in stochastic control theory is the sufficient condition for a minimum. The sufficient condition requires a suitably behaved solution $f$ of the dynamic programming equation (\ref{eqprg}) and a control law $\mathbf{u}^{\star}$ satisfying (\ref{verf3}). This result is called a \textit {verification theorem}.
\begin{theorem}(Sufficient conditions for optimality)
Let $f$ be a solution of (\ref{eqprg})-(\ref{boun1}) such that $f$ is in $C^2(G)\cap C(\overline{G}).$ Then:\\
(a) $f(\phi,\stackrel{\rightarrow}{x})\leq J(\phi,\stackrel{\rightarrow}{x},\mathbf{u})$ for any $\mathbf{u}\in\mathcal{U}$ and $(\phi,\stackrel{\rightarrow}{x})\in G.$\\
(b) If $\mathbf{u}^{\star}\in \mathcal{U},$  $J(\phi,\stackrel{\rightarrow}{x},\mathbf{u}^{\star})<\infty$ and
\begin{equation}\label{opticon}\label{verf3}
 \mathbf{A}^{\mathbf{u}^{\star}}f(\phi,\stackrel{\rightarrow}{x})+L^{\mathbf{u}^{\star}}(\phi,\stackrel{\rightarrow}{x})=\inf_{v\in U}\left[(\mathbf{A}^{v}f+L^v)(\phi,\stackrel{\rightarrow}{x})\right]
\end{equation}
for all $(\phi,\stackrel{\rightarrow}{x})\in G,$ then $f(\phi,\stackrel{\rightarrow}{x})=J(\phi,\stackrel{\rightarrow}{x},\mathbf{u}^{\star}).$ Thus $\mathbf{u}^{\star}$ is a optimal control, for all choices of initial data $(\phi,\stackrel{\rightarrow}{x})\in G.$
\end{theorem}
\proof (a). It is sufficient to consider those $\mathbf{u}$ for which $J(\phi,\stackrel{\rightarrow}{x},\mathbf{u})<\infty.$  The Chebishev inequality, (\ref{inflim}) and the boundedness of $\psi$ on $\partial^{\ast}(G),$ implies that $\mathbb{P}^{(\phi,\stackrel{\rightarrow}{x})}\left(\tau_{G}<\infty\right)=1.$ For each $v\in U,$ $(\phi,\stackrel{\rightarrow}{x})\in G,$
$0\leq(\mathbf{A}^{\mathbf{u}}f+L^{\mathbf{u}})(\phi,\stackrel{\rightarrow}{x}).$ We conclude the proof by using the Lemma (\ref{lemut}) replacing $M^{\mathbf{u}}$ by $L^{\mathbf{u}}.$\\
(b) The condition (\ref{condint}) imply that
$$
E^{(\phi,\stackrel{\rightarrow}{x})}\int_0^{\tau_{G}}|M^{\mathbf{u}}(S_t,S(t))|dt<\infty.
$$
For $\mathbf{u}=\mathbf{u}^{\star}$, we get $\mathbf{A}^{\mathbf{u}}f+L^{\mathbf{u}}(\phi,\stackrel{\rightarrow}{x})=0.$ Then, using Lemma (\ref{lemut})(b), we have $f(\phi,\stackrel{\rightarrow}{x})=J(\phi,\stackrel{\rightarrow}{x},\mathbf{u}^{\star}).$
\eop
\begin{definition}A point $(\phi,\stackrel{\rightarrow}{x})\in\partial(G)$ is called \textbf{regular} for $G$ (with respect to $(S_t,S(t))$) if
$$
\mathbb{P}^{(\phi,\stackrel{\rightarrow}{x})}\left(\tau_{G}=0\right)=1.
$$
Otherwise the point $(\phi,\stackrel{\rightarrow}{x})$ is called \textbf{irregular}.
\end{definition}
The verification theorem gives sufficient conditions for optimality. The following theorem gives necessary conditions for optimality, under sufficiently strong assumptions..
\begin{theorem}(The Hamilton-Jacobi-Bellman(HJB) equation)  Suppose that $\mathbb{P}^{(\phi,\stackrel{\rightarrow}{x})}\left(\tau_{G}<\infty\right)=1$ for each $(\phi,\stackrel{\rightarrow}{x})\in G.$
Define
$$
\Phi(\phi,\stackrel{\rightarrow}{x})=\inf_{\mathbf{u}}\left\{J^{\mathbf{u}}(\phi,\stackrel{\rightarrow}{x});\ \mathbf{u} \textmd{ Markov control}\right\}
$$
Suppose that $\Phi\in C^2\left(G\right)\bigcap C\left(\overline{G}\right)$ is bounded 
and that an optimal Markov control $\mathbf{u}^{\star}$ exists and that $\partial(G)$ is regular for $(S_t^{\mathbf{u}^{\star}},S^{\mathbf{u}^{\star}}(t)).$ Then
\begin{equation}\label{HJBe}
\inf_{v\in U}\left\{L^v(\phi,\stackrel{\rightarrow}{x})+\mathbf{A}^v\Phi(\phi,\stackrel{\rightarrow}{x})\right\}=0\ \forall\ (\phi,\stackrel{\rightarrow}{x})\in \left(G\right)
\end{equation}
and
\begin{equation}\label{bouncd}
\Phi(\phi,\stackrel{\rightarrow}{x})=\psi(\phi,\stackrel{\rightarrow}{x})\ \forall\ (\phi,\stackrel{\rightarrow}{x})\in\partial (G).
\end{equation}
The infimum in (\ref{HJBe}) is obtained if $v=\mathbf{u}^{\star}(\phi,\stackrel{\rightarrow}{x})$ where $\mathbf{u}^{\star}(\phi,\stackrel{\rightarrow}{x})$ is optimal. Equivalently
\begin{equation}\label{edp}
(L^{\mathbf{u}^{\star}(\phi,\stackrel{\rightarrow}{x})})(\phi,\stackrel{\rightarrow}{x})+(\mathbf{A}^{\mathbf{u}^{\star}(\phi,\stackrel{\rightarrow}{x})}\Phi)(\phi,\stackrel{\rightarrow}{x})=0\ \forall\ (\phi,\stackrel{\rightarrow}{x})\in G.
\end{equation}
\end{theorem}
\proof Since $\mathbf{u}^{\star}(\phi,\stackrel{\rightarrow}{x})$ is optimal, we obtain
\begin{equation}\begin{array}{rl}\label{optco}
\Phi(\phi,\stackrel{\rightarrow}{x})=&J^{\mathbf{u}^{\star}}{(\phi,\stackrel{\rightarrow}{x})}=E^{(\phi,\stackrel{\rightarrow}{x})}\left[\int_0^{\tau_{G}}L^{\mathbf{u}^{\star}(\phi,\stackrel{\rightarrow}{x})}(S_t,S(t))dt\right.+\\
&\!\!\!+\left.\psi(S_{\tau_{G}},S(\tau_{G}))\right].
\end{array}
\end{equation}
If $(\phi,\stackrel{\rightarrow}{x})\in\partial (G)$ then $\tau_{G}=0$ a.s. and we get (\ref{bouncd}). From (\ref{optco}) and Theorem \ref{possdrc} we obtain (\ref{edp}).\\
The proof is complete if we prove (\ref{HJBe}). Following \cite{O}, fix $(\phi,\stackrel{\rightarrow}{x})\in G$ and choose a Markov control $\mathbf{u}$. Let $\alpha\leq \tau_{G}$ be a bounded stopping time. 
Since
$$
J^{\mathbf{u}}(\phi,\stackrel{\rightarrow}{x})=E^{(\phi,\stackrel{\rightarrow}{x})}\left[\int_0^{\tau_{G}}L^{\mathbf{u}}(S_t,S(t))dt+\psi(S_{\tau_{G}},S(\tau_{G}))\right]
$$
using the Theorem \ref{strgmrkv} and the properties of the \textit{shift operator} $\theta_{\cdot}$ (see \cite{O} sections 7.2 and 9.3) we have
$$\begin{array}{rl}
E^{(\phi,\stackrel{\rightarrow}{x})}\left[J^{\mathbf{u}}(S_{\alpha},S(\alpha))\right]&\!\!\!\!=E^{(\phi,\stackrel{\rightarrow}{x})}\left[E^{(S_{\alpha},S(\alpha))}\left[\int_0^{\tau_{G}}L^{\mathbf{u}}(S_t,S(t))dt+\psi(S_{\tau_{G}},S(\tau_{G}))\right]\right]\\
=E^{(\phi,\stackrel{\rightarrow}{x})}&\!\!\!\!\left[E^{(\phi,\stackrel{\rightarrow}{x})}\left[\theta_{\alpha}\left(\int_0^{\tau_{G}}L^{\mathbf{u}}(S_t,S(t))dt+\psi(S_{\tau_{G}},S(\tau_{G}))\right)\left|\mathcal{F}_{\alpha}\right.\right]\right]\\
=E^{(\phi,\stackrel{\rightarrow}{x})}&\!\!\!\!\left[E^{(\phi,\stackrel{\rightarrow}{x})}\left[\theta_{\alpha}\left(\int_{\alpha}^{\tau_{G}}L^{\mathbf{u}}(S_t,S(t))dt+\psi(S_{\tau_{G}},S(\tau_{G}))\right)\left|\mathcal{F}_{\alpha}\right.\right]\right]\\
=E^{(\phi,\stackrel{\rightarrow}{x})}&\!\!\!\!\left[\int_0^{\tau_{G}}L^{\mathbf{u}}(S_t,S(t))dt+\psi(S_{\tau_{G}},S(\tau_{G}))-\int_0^{\alpha}L^{\mathbf{u}}(S_t,S(t))dt\right]\\
=J^{\mathbf{u}}(\phi,\stackrel{\rightarrow}{x})&\!\!\!\!-E^{(\phi,\stackrel{\rightarrow}{x})}\left[\int_0^{\alpha}L^{\mathbf{u}}(S_t,S(t))dt\right].
\end{array}
$$
Then
\begin{equation}\label{midcos}
J^{\mathbf{u}}(\phi,\stackrel{\rightarrow}{x})=E^{(\phi,\stackrel{\rightarrow}{x})}\left[\int_0^{\alpha}L^{\mathbf{u}}(S_t,S(t))dt\right]+E^{(\phi,\stackrel{\rightarrow}{x})}\left[J^{\mathbf{u}}(S_{\alpha},S(\alpha))\right].
\end{equation}
Now, we consider $W\subset G$ and $\alpha:=\inf\left\{t\geq 0; (S_t,S(t))\notin W\right\}.$ Suppose an optimal control $\mathbf{u}^{\star}(\phi,\stackrel{\rightarrow}{x})$ exists, let $v\in U$ arbitrary we define
$$
\mathbf{u}(\eta,\stackrel{\rightarrow}{y})=\left\{\begin{array}{rl}
v&\textmd{ if }(\eta,\stackrel{\rightarrow}{y})\in W,\\
\mathbf{u}^{\star}(\eta,\stackrel{\rightarrow}{y})&\textmd{ if }(\eta,\stackrel{\rightarrow}{y})\in G\setminus W.
\end{array}\right.
$$
Then
\begin{equation}
\Phi(S_{\alpha},S(\alpha))=J^{\mathbf{u}}(S_{\alpha},S(\alpha))=J^{\mathbf{u}^{\star}}(S_{\alpha},S(\alpha)),
\end{equation}
from this, (\ref{midcos}) and using the Dynkin formula we obtain
\begin{equation}\begin{array}{rl}
\Phi(\phi,\stackrel{\rightarrow}{x})&\leq J^{\mathbf{u}}(\phi,\stackrel{\rightarrow}{x})=E^{(\phi,\stackrel{\rightarrow}{x})}\left[\int_0^{\alpha}L^v(S_t,S(t))dt\right]+E^{(\phi,\stackrel{\rightarrow}{x})}\left[\Phi(S_{\alpha},S(\alpha))\right]\\
&=E^{(\phi,\stackrel{\rightarrow}{x})}\left[\int_0^{\alpha}L^v(S_t,S(t))dt\right]+\Phi(\phi,\stackrel{\rightarrow}{x})+\\
&+E^{(\phi,\stackrel{\rightarrow}{x})}\left[\int_0^{\alpha}\mathbf{A}^{v}\Phi(S_t,S(t))dt\right],
\end{array}
\end{equation}
therefore
$$
E^{(\phi,\stackrel{\rightarrow}{x})}\left[\int_0^{\alpha}\left(L^v(S_t,S(t))+\mathbf{A}^{v}\Phi(S_t,S(t))\right)dt\right]\geq0.
$$
Thus
$$
\frac{E^{(\phi,\stackrel{\rightarrow}{x})}\left[\int_0^{\alpha}\left(L^v(S_t,S(t))+\mathbf{A}^{v}\Phi(S_t,S(t))\right)dt\right]}{E^{(\phi,\stackrel{\rightarrow}{x})}[\alpha]}\geq0.
$$
Taking in account that $L^v(\cdot)$ and $\mathbf{A}^{v}(\cdot)$ are continuous, we obtain\\ $L^v(\phi,\stackrel{\rightarrow}{x})+\mathbf{A}^v(\phi,\stackrel{\rightarrow}{x})\geq0,$ as $W\downarrow(\phi,\stackrel{\rightarrow}{x}).$ From this and (\ref{edp}) we obtain (\ref{HJBe}).
\eop\\
\section{Example: An Optimal Portfolio Selection Pro- blem}
Let $S(t)$ denote the wealth of a person at time $t.$  The person has two investments. Let $P(t)$ be a \textit{risk free} investment:
$$
dP(t)=kP(t)dt.
$$
And the another investment is a \textit{risky} one:
$$
dS_1(t)=\mu S_1(t)dt+\sigma S_1(t)dW(t),
$$
and we assume that $k<\mu.$ At each instant $t$ the person can choose what fraction $u(t)$ of this wealth he will invest in the risky asset, then investing $1-u(t)$ in the risk free asset. Suppose that the past has influence over the wealth, $S(t),$ under the following SFDE
$$\begin{array}{rl}
dS(t)=&\mu u(t)\frac{S(t)}{1+\|S_t\|}dt+\sigma u(t)\frac{S(t)}{1+\|S_t\|}dW(t)+\\
&+k(1-u(t))\frac{S(t)}{1+\|S_t\|}dt=\\
&=\left(\mu u(t)+k(1-u(t))\right)\frac{S(t)}{1+\|S_t\|}dt+\sigma u(t)\frac{S(t)}{1+\|S_t\|}dW(t),
\end{array}
$$
and $(S_0,S(0))=(\phi,x)$ with $\|\phi\|>0$ and $x>0.$
By the Theorem \ref{lexuq1} there is a solution $S(t)$ with initial condition $(\phi,x).$\\
Assume that we do not allow any \textit{ borrowing} (i.e. require $u(t)\leq1$) and we do not allow any \textit{ shortselling} (i.e. require $0\leq u(t)$) and $\psi:[0,\infty)\rightarrow[0,\infty),$ $\psi(0)=0,$ (fixing this function) the problem is to find $\Xi(\phi,x)$ and a control $\mathbf{u}^{\star}=\mathbf{u}^{\star}(S_t,S(t)),$ $0\leq \mathbf{u}^{\star}\leq1,$ such that
$$
\Xi(\phi,x)=\sup\left\{J^{\mathbf{u}}(\phi,x):\ \mathbf{u}\textmd{ Markov control, }0\leq \mathbf{u}\leq1\right\}=J^{\mathbf{u}^{\star}}(\phi,x),
$$
where $J^{\mathbf{u}}(\phi,x)=E^{(\phi,x)}\left[\psi(S_{\tau_{G}}^{\mathbf{u}},S^{\mathbf{u}}(\tau_{G}))\right]$ and $\tau_{G}$ is the first exit time from $G=\left\{(\phi,x)\in V\times\mathbb{R}:\ x,\|\phi\|>0\textmd{ and } \frac{(\mu-k)^2}{2\sigma^2(1-p)}+\frac{k}{1+\|\phi\|}+\frac{\phi(0)^2-\phi(-r)^2}{p\|\phi\|^2}=0\right\}.$\\
We observe that
$$
\Xi=-\inf\left\{-J^{\mathbf{u}}(\phi,x)\right\}=-\inf\left\{E^{(\phi,x)}\left[-\psi(S_{\tau_{G}}^{\mathbf{u}},S^{\mathbf{u}}(\tau_{G}))\right]\right\},
$$
so $-\Xi$ coincides with the solution $\Phi$ of the problem (\ref{infprob}), but with $\psi$ replaced by $-\psi$ and $L=0.$ Thus, we see that the equation (\ref{eqprg}) for $\Phi$ gets for $\Xi$ the form
$$
\sup_v\left\{(\mathbf{A}^{v}f)(\phi,x)\right\}=0,\textmd{ for }(\phi,x)\in G;
$$
and
$$
f(\phi,x)=\psi(\phi,x)\textmd{ for } (\phi,x)\in\partial G.
$$
From (\ref{wkinfop}) the differential operator $\mathbf{A}^{v}$ has the form
$$\begin{array}{rl}
(\mathbf{A}^{v}f)(\phi,x)=&\!\!\!\!\frac{\partial f}{\partial x}(\phi,x)(\mu v+k(1-v))\frac{x}{1+\|\phi\|}+\frac{1}{2}\frac{\partial^2f}{\partial x^2}(\phi,x)\sigma^2v^2\frac{x^2}{(1+\|\phi\|)^2}+\\
&\!\!\!\!+\mathbf{\Gamma}f\left(\phi,x\right).
\end{array}
$$
Therefore, for each $(\phi,x)$ we try to find the value $v=(\phi,x)$ which maximizes the function
\begin{equation}\label{ex1hje1}\begin{array}{rl}
m(v)=&\!\!\!((\mu-k)v+k))\frac{x}{1+\|\phi\|}\frac{\partial f}{\partial x}(\phi,x)+\frac{1}{2}\sigma^2v^2\frac{x^2}{(1+\|\phi\|)^2}\frac{\partial^2f}{\partial x^2}(\phi,x)+\\
&\!\!\!+\mathbf{\Gamma}f\left(\phi,x\right).
\end{array}
\end{equation}
If $\frac{\partial f}{\partial x}>0$ and $\frac{\partial^2f}{\partial x^2}<0,$ the solution is
\begin{eqnarray}
v=\mathbf{u}(\phi,x)
&\!\!\!=&\!\!\!-\dfrac{(\mu-k)(1+\|\phi\|)\frac{\partial f}{\partial x}}{\sigma^2x\frac{\partial^2 f}{\partial x^2}}.\label{optimum1}
\end{eqnarray}
Replacing this in (\ref{ex1hje1}) we obtain the following boundary value problem
\begin{eqnarray}
&\!\!\!\!\!\!-&\!\!\!\!\frac{(\mu-k)^2}{2\sigma^2\frac{\partial^2 f}{\partial x^2}}(\frac{\partial f}{\partial x}(\phi,x))^2+k\frac{x}{1+\|\phi\|}\frac{\partial f}{\partial x}(\phi,x)+\nonumber\\
&\!\!\!\!\!\!\!+&\!\!\!\!\!\!\mathbf{\Gamma}f\left(\phi,x\right)=0\label{valueproeq1}\\
f&\!\!\!\!\!\!(&\!\!\!\!\!\phi,x)=\psi(\phi,x)\textmd{ for } (\phi,x)\in\partial G\label{valueprobond1}
\end{eqnarray}
Consider $\psi(\phi,x)=x^p$ where $0<p<1.$\\
We try to find a solution of (\ref{valueproeq1}) and (\ref{valueprobond1}) of the form
$$f(\phi,x)=\|\phi\|^2x^p.$$
Substituting into (\ref{valueproeq1}) and using the definition of $\mathbf{\Gamma}$ we obtain $\frac{p(\mu-k)^2\|\phi\|^2}{2\sigma^2(1-p)}+\frac{kp\|\phi\|^2}{1+\|\phi\|}+\phi(0)^2-\phi(-r)^2=0.$\\
Using (\ref{optimum1}) we obtain the optimal control
$$
\mathbf{u}^{\star}(\phi,x)=\frac{(\mu-k)(1+\|\phi\|)}{\sigma^2(1-p)}.
$$
If $0<\frac{(\mu-k)(1+\|\phi\|)}{\sigma^2(1-p)}<1,$ this $\mathbf{u}^{\star}$ is the solution to the problem.


\begin{thebibliography}{AAA}
\bibitem{SE2}{}Arriojas, M.,  Hu, Y.,  Mohammed, S.-E. A. and  Pap, G., A Delayed Black and Scholes Formula, \textit{Stochastic Analysis and Applications}, \textbf{25}:2  (2007), 471-492.
\bibitem{D}{}Dynkin, E. B., \textit{Markov Process,} Vol I, Die Grundlehreu der Math. Wissenschaften, Springer-Verlag, 1965.
\bibitem{DZ}{}Da Prato, G. and Zabczyk, J., \textit{Stochastics Equations in Infinite Dimensions}, Cambridge University Press, 1992.
\bibitem{FR}{}Fleming, W. H. and Rishel, R. W.,\textit{Deterministic and Stochastic Control,} Springer-Verlag, 1975.
\bibitem{GS}{}Gihman, I. I. and Skorohod, A. I., {\it Stochastic Differential Equations}, Springer-Verlag, 1972.
\bibitem{HO}{}Hu, Y. and {\O}ksendal, B., Fractional white noise calculus and applications to finance. {\it Infinite Dimensional Analysis, Quantum Probability and Related Topics}, \textbf{6}:1 (2003), 1-32.
\bibitem{IS}{}Ivanov,  A. F. and Swishchuk, A. V., Optimal Control of Stochastic Differential Delay Equations with Applications in Economics, \textit{International Journal of Qualitative Theory of Differential Equations and Applications}
\textbf{2}:2 (2008), 201-213
\bibitem{KS}{} Karatzas, I. and Shreve, S. E., {\it Brownian Motion and Stochastic Calculus}, Second Edition, Springer, NY, 1991.
\bibitem{K}{}Kushner, H. J., On the Stability of Process Defined by Stochastic Difference-Differential Equations, {\it Journal of Differential Equations}, \textbf{4}, (1968), 424-443.
\bibitem{SE3}{}Mohammed, S.-E. A., \textit{Stochastic Functional Differential Equations.} Research Notes in Mathematics No. 99, Pitman Books Ltd., London, 1984.
\bibitem{O}{} {\O}ksendal, B., {\it Stochastic Differential Equations. An Introduction with Applications} Springer-Verlag, Sixth. Ed. 2003.
\bibitem{R}{}Ramsey, F. P., A mathematical theory of savings, {\it The Economic Journal}, \textbf{38}:152 (1928), 543-549,.
\bibitem{SK}{}Schoenmakers, J., and Kloeden, P., Robust option replication for a Black- Scholes model extended with nondeterministic trends. {\it Journal of Applied Mathematics and Stochastic Analysis}, \textbf{12}:2 (1999), 113-120.
\bibitem{YSE}{}Yan, F., and Mohammed, S.-E. A., A Stochastic Calculus for Systems with Memory, {\it Stochastic Analysis and Applications}, \textbf{23}:3 (2005) 613-657.
\end{thebibliography}
\end{document}